\documentclass[twoside,reqno]{amsart}
\usepackage{latexsym}%%triangles:\lhd\rhd\unlhd\unrhd%
\newcommand{\qdn}{\hspace*{-1.5mm}}

\newcommand{\xqdn}{\hspace*{-5.0mm}}
\newcommand{\xxqdn}{\hspace*{-10mm}}

\newcommand{\fns}{\footnotesize}

\newcommand{\ffnk}[4]{\left[\qdn\ba{#1}#3\\[0.7mm]#4\ea{\!\Big|\:#2}\right]}

\newcommand{\nnm}{\nonumber}
\newcommand{\be}{\begin{equation}}
\newcommand{\ee}{\end{equation}}
\newcommand{\ba}{\begin{array}}
\newcommand{\ea}{\end{array}}
\newcommand{\bmn}{\begin{eqnarray}}
\newcommand{\emn}{\end{eqnarray}}
\newcommand{\bnm}{\begin{eqnarray*}}
\newcommand{\enm}{\end{eqnarray*}}
\newcommand{\bln}{\begin{subequations}}
\newcommand{\eln}{\end{subequations}}

\newtheorem{thm}{Theorem}%[section]

\newtheorem{exam}{Example}

\newtheorem{entry}{Entry}%%%%%%%%%%%%%%%%

\newcommand{\bbtm}[4]{\bibitem{kn:#1}{#2,}~{#3,}~{#4.}}
\newcommand{\cito}[1]{\cite{kn:#1}}
\newcommand{\citu}[2]{\cite[#2]{kn:#1}}

\begin{document} %%%%%%%%%% This paper is published in %%%%%%%
{\fns% \today\hfill\copyright%% Printed in China} %%%%%%%%%%%%%%%
%%%%%%%%%%%%%%%%%%%%%%%%%%%%%%%%%%%%%%%%%%%%%%%%%%%%%%%%%%%%%%
\title{Extensions of Ramanujan's two formulas for $1/\pi$}
\author{$^a$Chuanan Wei, $^b$Dianxuan Gong}
\dedicatory{
$^A$Department of Information Technology\\
  Hainan Medical College,  Haikou 571101, China\\
         $^B$College of Sciences\\
             Hebei Polytechnic University, Tangshan 063009, China}
\thanks{\emph{Email addresses}:
 weichuanan@yahoo.com.cn(C. Wei),  gongdianxuan@yahoo.com.cn(D.Gong)}

\address{ }
\footnote{\emph{2010 Mathematics Subject Classification}: Primary
33C20 and Secondary 40A15, 65B10}

\keywords{Hypergeometric series; Ramanujan's two formulas for
$1/\pi$; Summation formula for $1/\pi$ with free parameters}

\begin{abstract}
In terms of the hypergeometric method, we establish the extensions
of two formulas for $1/\pi$ due to Ramanujan \cito{ramanujan}.
Further, other five summation formulas for $1/\pi$ with free
parameters are also derived in the same way.
\end{abstract}

%%%%%%%%%%%%%%%%%%%%%%%%%%%%%%%%%%%%%%%%%%%%%%%%%%%%%%%%%%%%%%%%%%%
\maketitle\thispagestyle{empty}%%%%%%%%%%%%%%%%%%%%%%%%%%%%%%%%%%%%
\markboth{Chuanan Wei, Dianxuan Gong}%%%%%%%%%%%%%%%%%%%%%%%%%%%%
         {Extensions of Ramanujan's two formulas for $1/\pi$}

%%%%%%%%%%%%%%%%%%%%%%%%%%%%%%%%%%%%%%%%%%%%%%%%%%%%%%%%%%%%%%%%%%%
%%%%%%%%%%%%%%%%%%%%%%%%%%%%%%%%%%%%%%%%%%%%%%%%%%%%%%%%%%%%%%%%%%%
%%%%%%%%%%%%%%%%%%%%%%%%%%%%%%%%%%%%%%%%%%%%%%%%%%%%%%%%%%%%%%%%%%%
\section{Introduction}

 For a complex number $x$ and an integer $n$, define the shifted
factorial by
 \[(x)_n=\Gamma(x+n)/\Gamma(x)\]
where $\Gamma$-function is well-defined:
\[\Gamma(x)=\int_{0}^{\infty}t^{x-1}e^{-t}dt\:\:\text{with}\:\:Re(x)>0.\]
For centuries, the study of $\pi$-formulas attracts many
mathematicians. The corresponding results can be found in
\cito{adanmchik}-\cito{bailey-b}, \cito{baruch-a}-\cito{chu-c} and
\cito{glaisher}-\cito{zudilin-b}. Thereinto, two formulas for
$1/\pi$ due to Ramanujan \cito{ramanujan} can be stated as
 \bmn
 &&\sum_{k=0}^{\infty}\frac{(\frac{1}{2})_k^3}{(k!)^3}\frac{6k+1}{4^k}
=\frac{4}{\pi}, \label{ram-a}\\
 &&\sum_{k=0}^{\infty}\frac{(\frac{1}{2})_k(\frac{1}{4})_k(\frac{3}{4})_k}{(k!)^3}\frac{8k+1}{9^k}
=\frac{2\sqrt{3}}{\pi}. \label{ram-b}
 \emn

Following  Bailey~\cito{bailey}, define the hypergeometric series
 by
\[_{1+r}F_s\ffnk{cccc}{z}{a_0,&a_1,&\cdots,&a_r}
{&b_1,&\cdots,&b_s} \:=\:\sum_{k=0}^\infty
\frac{(a_{0})_{k}(a_{1})_{k}\cdots(a_{r})_{k}}
 {k!(b_{1})_{k}\cdots(b_{s})_{k}}z^k.\]
Then the identity due to Gessel-Stanton \citu{gessel}{Eq. (1.7)} and
Gasper's identity (cf. \citu{gasper}{Eq. (5.23)}) can be expressed
as
 \bmn
&&_7F_6\ffnk{cccccccc}{1}{a,1+\frac{a}{3},b,1-b,c,\frac{1}{2}+a-c+n,-n}
 {\frac{a}{3},\frac{2+a-b}{2},\frac{1+a+b}{2},1+a-2c,1+a+2n,2c-a-2n}\nnm\\
 &&\:\,=\:\,\frac{(\frac{1+a}{2})_n(1+\frac{a}{2})_n(\frac{1+a+b}{2}-c)_n(1+\frac{a-b}{2}-c)_n}
 {(\frac{1+a+b}{2})_n(1+\frac{a-b}{2})_n(\frac{1+a}{2}-c)_n(1+\frac{a}{2}-c)_n},
 \label{hypergeotric-a}\\
&&_7F_6\ffnk{cccccccc}{1}{3a,1+\frac{3a}{4},\frac{1-3b}{2},\frac{2-3b}{2},3b,2a+b+n,-n}
 {\frac{3a}{4},\frac{1+3a+3b}{2},\frac{3a+3b}{2},1+a-b,1-3a-3b-3n,1+3a+3n}\nnm\\
 &&\:\,=\:\,\frac{(a+2b)_n(a+\frac{1}{3})_n(a+\frac{2}{3})_n(a+1)_n}
 {(1+a-b)_n(a+b)_n(a+b+\frac{1}{3})_n(a+b+\frac{2}{3})_n}.
 \label{hypergeotric-b}
 \emn

The main aim of the paper is to explore the relations of
hypergeometric series and $\pi$-formulas. Four summation formulas
for $1/\pi$ with free parameters including the extension of
\eqref{ram-a} will be derived from \eqref{hypergeotric-a} in section
2. Three summation formulas for $1/\pi$ with free parameters
including the extension of \eqref{ram-b} will be deduced from
\eqref{hypergeotric-b} in section 3.

\section{Summation formulas for $1/\pi$ with free parameters\\ implied by the Gessel-Stanton identity}

Letting $n\to\infty$ for \eqref{hypergeotric-a}, we obtain the
following equation:
 \bmn
&&_5F_4\ffnk{cccccccc}{\frac{1}{4}}{a,1+\frac{a}{3},b,1-b,c}
 {\frac{a}{3},\frac{2+a-b}{2},\frac{1+a+b}{2},1+a-2c}
  \nnm\\\label{hypergeometric-aa}&&\:\,=\:
 \frac{\Gamma(\frac{1+a+b}{2})\Gamma(1+\frac{a-b}{2})\Gamma(\frac{1+a}{2}-c)\Gamma(1+\frac{a}{2}-c)}
 {\Gamma(\frac{1+a}{2})\Gamma(1+\frac{a}{2})\Gamma(\frac{1+a+b}{2}-c)\Gamma(1+\frac{a-b}{2}-c)}.
 \emn

Choosing $a=\frac{1}{2}+2p$, $b=\frac{1}{2}+2q$ and
$c=\frac{1}{2}+r$ in \eqref{hypergeometric-aa}, we achieve the
extension of \eqref{ram-a}.

\begin{thm}\label{thm-a} For $p,q,r\in \mathbb{Z}$ with $\min\{p+q,p-q\}\geq0$,
 there holds the summation formula for $1/\pi$ with free parameters:
\bnm
 \quad\frac{1}{\pi}&&\xqdn=\:\frac{(\frac{1}{2})_{p+q-r}(\frac{1}{2})_{p-q-r}}
{(\frac{1}{2})_{r}}\\
&&\xqdn\times\:\sum_{k=0}^{\infty}\frac{(\frac{1}{2})_{k+2p}(\frac{1}{2})_{k+2q}(\frac{1}{2})_{k-2q}(\frac{1}{2})_{k+r}}
 {k!(k+p+q)!(k+p-q)!(\frac{1}{2})_{k+2p-2r}}\frac{6k+4p+1}{4^{k+r+1}}.
 \enm
\end{thm}

When $p=q=r=0$, Theorem \ref{thm-a} reduces to Ramanujan's formula
for $1/\pi$ given by \eqref{ram-a} exactly. Other two examples of
the same type are displayed as follows.

\begin{exam}[$p=q=1$, $r=2$ in Theorem \ref{thm-a}]
 \bnm
\frac{256}{3\pi}=\sum_{k=0}^{\infty}\frac{(\frac{5}{2})_k^3}{(k!)^2(k+2)!}\frac{6k+5}{4^k}.
 \enm
\end{exam}

\begin{exam}[$p=q=2$, $r=4$ in Theorem \ref{thm-a}]
 \bnm\quad\:
\frac{16384}{315\pi}=\sum_{k=0}^{\infty}\frac{(\frac{9}{2})_k^3}{(k!)^2(k+4)!}\frac{2k+3}{4^k}.
 \enm
\end{exam}

Making $a=\frac{3}{2}+2p$, $b=\frac{3}{2}+2q$ and $c=\frac{1}{2}+r$
in \eqref{hypergeometric-aa}, we attain the identity.

\begin{thm}\label{thm-b} For $p,q,r\in \mathbb{Z}$ with $\min\{p+q+1,p-q\}\geq0$,
 there holds the summation formula for $1/\pi$ with free parameters:
\bnm
 \quad\frac{1}{\pi}&&\xqdn=\:\frac{(\frac{1}{2})_{p+q-r+1}(\frac{1}{2})_{p-q-r}}{(\frac{1}{2})_{r}}\\
&&\xqdn\times\:\sum_{k=0}^{\infty}\frac{(\frac{3}{2})_{k+2p}(\frac{3}{2})_{k+2q}(-\frac{1}{2})_{k-2q}(\frac{1}{2})_{k+r}}
 {k!(k+p+q+1)!(k+p-q)!(\frac{3}{2})_{k+2p-2r}}\frac{6k+4p+3}{4^{k+r+1}}.
 \enm
\end{thm}

Two examples from Theorem \ref{thm-b} are laid out as follows.

\begin{exam}[$p=q=0$, $r=1$ in Theorem \ref{thm-b}]
 \bnm
\frac{32}{3\pi}=\sum_{k=0}^{\infty}\frac{(\frac{3}{2})_k^3}{(k!)^2(k+1)!}\frac{2k+1}{4^k}.
 \enm
\end{exam}

\begin{exam}[$p=q=1$, $r=3$ in Theorem \ref{thm-b}]
 \bnm\quad\:
\frac{2048}{15\pi}=\sum_{k=0}^{\infty}\frac{(\frac{7}{2})_k^3}{(k!)^2(k+3)!}\frac{6k+7}{4^k}.
 \enm
\end{exam}

Taking $a=\frac{1}{2}+2p$, $b=\frac{1}{2}+2q$ in
\eqref{hypergeometric-aa} and then letting $c\to-\infty$, we get the
identity.

\begin{thm}\label{thm-c} For $p,q\in \mathbb{Z}$ with $\min\{p+q,p-q\}\geq0$,
 there holds the summation formula for $1/\pi$ with free parameters:
\bnm
 \quad\frac{4^{p+1}}{\sqrt{2}\,\pi}=
\sum_{k=0}^{\infty}\frac{(\frac{1}{2})_{k+2p}(\frac{1}{2})_{k+2q}(\frac{1}{2})_{k-2q}}
 {k!(k+p+q)!(k+p-q)!}\frac{6k+4p+1}{(-8)^k}.
 \enm
\end{thm}

Two examples from Theorem \ref{thm-c} are displayed as follows.

\begin{exam}[$p=q=0$ in Theorem \ref{thm-c}]
 \bnm \xxqdn
\frac{2\sqrt{2}}{\pi}=\sum_{k=0}^{\infty}\frac{(\frac{1}{2})_k^3}{(k!)^3}\frac{6k+1}{(-8)^k}.
 \enm
\end{exam}

\begin{exam}[$p=q=1$ in Theorem \ref{thm-c}]
 \bnm \quad
  \frac{32\sqrt{2}}{3\pi}=\sum_{k=0}^{\infty}\frac{(\frac{5}{2})_k^2(-\frac{3}{2})_k}{(k!)^2(k+2)!}\frac{6k+5}{(-8)^k}.
 \enm
\end{exam}

Setting $a=\frac{3}{2}+2p$, $b=\frac{3}{2}+2q$ in
\eqref{hypergeometric-aa} and then letting $c\to-\infty$, we gain
the identity.

\begin{thm}\label{thm-d} For $p,q\in \mathbb{Z}$ with $\min\{p+q+1,p-q\}\geq0$,
 there holds the summation formula for $1/\pi$ with free parameters:
\bnm
 \quad\frac{4^{p+2}}{\sqrt{2}\,\pi}=
\sum_{k=0}^{\infty}\frac{(\frac{3}{2})_{k+2p}(\frac{3}{2})_{k+2q}(-\frac{1}{2})_{k-2q}}
 {k!(k+p+q+1)!(k+p-q)!}\frac{6k+4p+3}{(-8)^k}.
 \enm
\end{thm}

Two examples from Theorem \ref{thm-d} are laid out as follows.

\begin{exam}[$p=q=0$ in Theorem \ref{thm-d}]
 \bnm
 \frac{8\sqrt{2}}{3\pi}=\sum_{k=0}^{\infty}\frac{(\frac{3}{2})_k^2(-\frac{1}{2})_k}{(k!)^2(k+1)!}\frac{2k+1}{(-8)^k}.
 \enm
\end{exam}

\begin{exam}[$p=q=1$ in Theorem \ref{thm-d}]
 \bnm \quad
  \frac{128\sqrt{2}}{15\pi}=\sum_{k=0}^{\infty}\frac{(\frac{7}{2})_k^2(-\frac{5}{2})_k}{(k!)^2(k+3)!}\frac{6k+7}{(-8)^k}.
 \enm
\end{exam}

\section{Summation formulas for $1/\pi$ with free parameters\\ implied by Gasper's identity}

Letting $n\to\infty$ for \eqref{hypergeotric-b}, we obtain the
following equation:
 \bmn
&&_5F_4\ffnk{cccccccc}{\frac{1}{9}}{3a,1+\frac{3a}{4},3b,\frac{1-3b}{2},\frac{2-3b}{2}}
 {\frac{3a}{4},1+a-b,\frac{1+3a+3b}{2},\frac{3a+3b}{2}}
  \nnm\\\label{hypergeometric-bb}&&\:\,=\:
 \frac{\Gamma(1+a-b)\Gamma(a+b)\Gamma(a+b+\frac{1}{3})\Gamma(a+b+\frac{2}{3})}
 {\Gamma(a+2b)\Gamma(a+\frac{1}{3})\Gamma(a+\frac{2}{3})\Gamma(a+1)}.
 \emn

Choosing $a=\frac{1}{6}+p$ and $b=\frac{1}{6}+q$ in
\eqref{hypergeometric-bb}, we achieve the extension of
\eqref{ram-b}.

\begin{thm}\label{thm-e} For $p,q\in \mathbb{Z}$ with $\min\{p+q,p-q\}\geq0$,
 there holds the summation formula for $1/\pi$ with free parameters:
\bnm
 \quad\frac{2(-1)^q}{3^{3q-\frac{1}{2}}\pi}=\Big(\frac{1}{2}\Big)_{p+2q}
\sum_{k=0}^{\infty}\frac{(\frac{1}{2})_{k+3p}(\frac{1}{2})_{k+3q}(\frac{1}{2})_{2k-3q}}
 {k!(k+p-q)!(2k+3p+3q)!}\frac{8k+6p+1}{9^k}.
 \enm
\end{thm}

When $p=q=0$, Theorem \ref{thm-e} reduces to Ramanujan's formula for
$1/\pi$ offered by \eqref{ram-b} exactly. Other two examples of the
same type are displayed as follows.

\begin{exam}[$p=q=1$ in Theorem \ref{thm-e}]
 \bnm
\frac{1024\sqrt{3}}{405\pi}=\sum_{k=0}^{\infty}\frac{(\frac{7}{2})_k(-\frac{5}{4})_k(-\frac{3}{4})_k}
{(k!)^2(k+3)!}\frac{8k+7}{9^k}.
 \enm
\end{exam}

\begin{exam}[$p=q=2$ in Theorem \ref{thm-g}]
\bnm \qquad\:
\frac{524288\sqrt{3}}{7577955\pi}=\sum_{k=0}^{\infty}\frac{(\frac{13}{2})_k(-\frac{11}{4})_k(-\frac{9}{4})_k}
{(k!)^2(k+6)!}\frac{8k+13}{9^k}.
 \enm
\end{exam}

Taking $a=\frac{1}{2}+p$ and $b=\frac{1}{2}+q$ in
\eqref{hypergeometric-bb}, we attain the identity.

\begin{thm}\label{thm-f} For $p,q\in \mathbb{Z}$ with $\min\{p+q,p-q\}\geq0$,
 there holds the summation formula for $1/\pi$ with free parameters:
\bnm
 \qquad\frac{4(-1)^q}{3^{3q+\frac{1}{2}}\pi}=\Big(\frac{1}{2}\Big)_{p+2q+1}
\sum_{k=0}^{\infty}\frac{(\frac{3}{2})_{k+3p}(\frac{3}{2})_{k+3q}(-\frac{1}{2})_{2k-3q}}
 {k!(k+p-q)!(2k+3p+3q+2)!}\frac{8k+6p+3}{9^k}.
 \enm
\end{thm}

Two examples from Theorem \ref{thm-f} are laid out as follows.

\begin{exam}[$p=q=0$ in Theorem \ref{thm-f}]
 \bnm
\frac{16\sqrt{3}}{3\pi}=\sum_{k=0}^{\infty}\frac{(\frac{3}{2})_k(-\frac{1}{4})_k(\frac{1}{4})_k}
{(k!)^2(k+1)!}\frac{8k+3}{9^k}.
 \enm
\end{exam}

\begin{exam}[$p=q=1$ in Theorem \ref{thm-f}]
\bnm \qquad\:
\frac{8192\sqrt{3}}{8505\pi}=\sum_{k=0}^{\infty}\frac{(\frac{9}{2})_k(-\frac{7}{4})_k(-\frac{5}{4})_k}
{(k!)^2(k+4)!}\frac{8k+9}{9^k}.
 \enm
\end{exam}

Setting $a=\frac{5}{6}+p$ and $b=\frac{5}{6}+q$ in
\eqref{hypergeometric-bb}, we get the identity.

\begin{thm}\label{thm-g} For $p,q\in \mathbb{Z}$ with $\min\{p+q+1,p-q\}\geq0$,
 there holds the summation formula for $1/\pi$ with free parameters:
\bnm
 \qquad\frac{8(-1)^q}{3^{3q+\frac{5}{2}}\pi}=\Big(\frac{1}{2}\Big)_{p+2q+2}
\sum_{k=0}^{\infty}\frac{(\frac{5}{2})_{k+3p}(\frac{5}{2})_{k+3q}(-\frac{3}{2})_{2k-3q}}
 {k!(k+p-q)!(2k+3p+3q+4)!}\frac{8k+6p+5}{9^k}.
 \enm
\end{thm}

Two examples from Theorem \ref{thm-g} are displayed as follows.

\begin{exam}[$p=q=0$ in Theorem \ref{thm-g}]
 \bnm
\frac{128\sqrt{3}}{27\pi}=\sum_{k=0}^{\infty}\frac{(\frac{5}{2})_k(-\frac{3}{4})_k(-\frac{1}{4})_k}
{(k!)^2(k+2)!}\frac{8k+5}{9^k}.
 \enm
\end{exam}

\begin{exam}[$p=q=1$ in Theorem \ref{thm-g}]
\bnm \qquad\:
\frac{65536\sqrt{3}}{229635\pi}=\sum_{k=0}^{\infty}\frac{(\frac{11}{2})_k(-\frac{9}{4})_k(-\frac{7}{4})_k}
{(k!)^2(k+5)!}\frac{8k+11}{9^k}.
 \enm
\end{exam}

\textbf{Remark:} With the change of the parameters, Theorems
\ref{thm-a}-\ref{thm-g} can produce more concrete formulas for
$1/\pi$. We shall not lay them out here.
%%%%%%%%%%%%%%%%%%%%%%%%%%%%%%%%%%%%%%%%%%%%%%%%%%%%%%%%%%%%%%%%%%%

%%%%%%%%%%%%%%%%%%%%%%%%%%%%%%%%%%%%%%%%%%%%%%%%%%%%%%%%%%%%%%%%%%%
%%%%%%%%%%%%%%%%%%%%%%%%%%%%%%%%%%%%%%%%%%%%%%%%%%%%%%%%%%%%%%%%%%%
%%%%%%%%%%%%%%%%%%%%%%%%%%%%%%%%%%%%%%%%%%%%%%%%%%%%%%%%%%%%%%%%%%%

\end{document}